\documentclass[11pt,reqno]{amsart}
\usepackage{amscd,amssymb,amsmath,amsthm}
\usepackage[arrow,matrix]{xy}
\usepackage{graphicx}
\usepackage{cite}
\usepackage{geometry}
\newtheorem{thm}[subsection]{Theorem}
\newtheorem{lemma}[subsection]{Lemma}
\newtheorem{pro}[subsection]{Proposition}

\numberwithin{equation}{section} 

\newcommand{\G}{{\mathcal G}}

\newcommand{\s}{{\sigma}}

\def \G {\Gamma}
\def \s {\sigma}

\newcommand{\bea}{\begin{eqnarray}}
\newcommand{\eea}{\end{eqnarray}}

\begin{document}
\title[Lyapunov operator and Gibbs measure]{Lyapunov operator $\mathcal L$ with degenerate kernel and Gibbs measures}
\author{Yu.\ Kh.\ Eshkabilov, F. H. Haydarov }

\address{Yu.\ Kh.\ Eshkabilov\\ Karshi State University, Kashkadarya, Uzbekistan}
\email {yusup62@mail.ru}

\address{F.\ H.\ Haydarov\\ National University of Uzbekistan,
Tashkent, Uzbekistan.} \email {haydarov\_imc@mail.ru}

\begin{abstract} In this paper we'll give a connection between four competing
interactions (external field, nearest neighbor, second neighbors
and triples of neighbors) of models with uncountable (i.e.
$[0,1]$) set of spin values on the Cayley tree of order two and
Lyapunov integral equation. Also we'll study fixed points of
Lyapunov operator with degenerate kernel which each fixed point of
the operator is correspond to a {\it translation-invariant} Gibbs
measure.
\end{abstract}
\maketitle



\section{Preliminaries}
Spin systems on lattices are a large class of systems considered
in statistical mechanics. Some of them have a real physical
meaning, others are studied as suitably simplified models of more
complicated systems \cite{newg}, \cite{11}.

The various partial cases of Ising model have been investigated in
numerous works, for example, the case $J_{3}=\alpha=0$ was
considered in \cite{gr15} and \cite{gr16}, the exact solutions of
an Ising model with competing restricted interactions with zero
external field was presented. In \cite{gr17} it is proved that
there are two \emph{translation-invariant} and uncountable number
of distinct \emph{non-translation-invariant} extreme Gibbs
measures. In \cite{gr9} the phase transition problem was solved
for $\alpha=0,\ J\cdot J_{1}\cdot J_{3}\neq 0$ and for $J_{3}=0,\
\alpha\cdot J\cdot J_{1}\neq 0$ as well. In \cite{p1a} it's
considered Ising model with four competing interactions (i.e.,
$J\cdot J_{1}\cdot J_{3}\cdot \alpha\neq 0$ ) on the Cayley tree
of order two. Mainly these papers are devoted to models with a
\emph{\textbf{finite}} set of spin values and in \cite{uar} given
other important results on a Cayley tree. In \cite{6} the Potts
model with a {\it \textbf{countable}} set of spin values on a
Cayley tree is considered and it was showed that the set of
translation-invariant splitting Gibbs measures of the model
contains at most one point, independently on parameters of the
Potts model with countable set of spin values on the Cayley tree.

It has been considering Gibbs measures for models with
\emph{\textbf{uncountable}} set of spin values for last five
years. Until now it has been considered models with
nearest-neighbor interactions $(J_{3}=J=\alpha=0,\ J_{1}\neq 0)$
and with the set $[0,1]$ of spin values on a Cayley tree and
gotten following results: "Splitting Gibbs measures" of the model
on a Cayley tree of order $k$ is described by solutions of a
nonlinear integral equation.
 For $k=1$ it's shown that the integral equation has a
unique solution (i.e., there is a unique Gibbs measure).
 For periodic splitting Gibbs measures it was found a sufficient condition under which
 the measure is unique and was proved existence of phase
 transitions on a Cayley tree of order $k\geq 2$ (see
 \cite{ehr2013}, \cite{ehr2012},\cite{rh2015}).

  In \cite{newh} it's considered splitting Gibbs measures for four competing
interactions i.e. ($J\cdot J_{1}\cdot J_{3}\cdot \alpha\neq 0$) of
models with uncountable set of spin values on the Cayley tree of
order two. In this paper we'll give a connection between Gibbs
measures for a given model and solutions of Lyupanov integral
equations. Also we'll study fixed points of Lyapunov operator with
degenerate kernel. Each fixed point of the operator is correspond
to a {\it translation-invariant} Gibbs measure.

 A Cayley tree $\G^k=(V,L)$ of order $k\in
\mathbb{N}$ is an infinite homogeneous tree, i.e., a graph without
cycles, with exactly $k+1$ edges incident to each vertices. Here
$V$ is the set of vertices and $L$ that of edges (arcs). Two
vertices $x$ and $y$ are called nearest neighbors if there exists
an edge $l\in L$ connecting them. We will use the notation
$l=\langle x,y\rangle$. The distance $d(x,y), x,y \in V$ on the
Cayley tree is defined by the formula

$$d(x,y)=\min\{d |\ x=x_{0},x_{1},...,x_{d-1},x_{d}=y\in V \ \emph{such that the pairs}$$
$$<x_{0},x_{1}>,...,<x_{d-1},x_{d}> \emph{are neighboring vertices}\}.$$\vskip
0.3 truecm

Let $x^{0}\in V$ be a fixed and we set

$$W_{n}=\{x\in V\ |\ d(x,x^{0})=n\}, \,\,\,\,\ V_{n}=\{x\in V\ |\ d(x,x^{0})\leq n\},$$

$$L_{n}=\{l=<x,y>\in L\ |\ x,y \in V_{n}\},$$\vskip
0.3 truecm
 The set of the direct successors of $x$ is denoted by $S(x),$
 i.e.
 $$S(x)=\{y\in W_{n+1}|\ d(x,y)=1\}, \ x\in W_{n}.$$
 We observe that for any vertex $x\neq x^{0},\ x$ has $k$ direct
 successors and $x^{0}$ has $k+1$. The vertices $x$ and $y$ are called second neighbor which is denoted by
 $>x,y<,$ if there exist a vertex $z\in V$ such that
 $x$, $z$ and $y$, $z$ are nearest neighbors. We will consider only second neighbors $> x, y <,$ for which there
exist $n$ such that $x, y \in W_n$. Three vertices $x,\ y$ and $z$
are called a triple of neighbors and they are denoted by $< x, y,
z>,$ if $< x, y >,\ < y, z >$ are nearest neighbors and $x,\ z \in
W_n,\ y \in W_{n-1}$, for some $n \in \mathbb{N}$.

Now we consider models with four competing interactions where the
spin takes values in the set
 $[0,1]$. For some set $A\subset V$ an arbitrary function $\s_A:A\to
[0,1]$ is called a configuration and the set of all configurations
on $A$ we denote by $\Omega_A=[0,1]^A$. Let $\sigma(\cdot)$ belong
to $\Omega_{V}=\Omega$ and $\xi_{1}:(t,u,v)\in[0,1]^{3}\to
\xi_{1}(t,u,v)\in R$,  $\xi_{i}: (u,v)\in [0,1]^2\to
\xi_{i}(u,v)\in R, \ i\in \{2,3\}$ are given bounded, measurable
functions. Then we consider the model with four competing
interactions on the Cayley tree which is defined by following
Hamiltonian

$$H(\sigma)=-J_{3}\sum_{<x,y,z>}\xi_{1}\left(\sigma(x),\sigma(y),\sigma(z)\right)
-J\sum_{>x,y<}\xi_{2}\left(\sigma(x),\sigma(z)\right)$$
\begin{equation}\label{e1}
-J_{1}\sum_{<x,y>}\xi_{3}\left(\sigma(x),\sigma(y)\right)-\alpha\sum_{x\in
V}\sigma(x),
\end{equation}\\
 where the sum in the first term ranges all triples of
neighbors, the second sum ranges all second neighbors, the third
sum ranges all nearest neighbors and  $J, J_{1}, J_{3},\alpha\in
R\setminus \{0\}$.
 Let $h: [0,1]\times V\setminus \{x^{0}\}\rightarrow \mathbb{R}$ and
  $|h(t,x)|=|h_{t,x}|<C$ where $x_{0}$ is a root of Cayley tree and $C$ is a
constant which does not depend on $t$.  For some $n\in\mathbb{N},$
$\sigma_n:x\in V_n\mapsto \sigma(x)$ and $Z_n$ is the
corresponding partition function we consider the probability
distribution $\mu^{(n)}$ on $\Omega_{V_n}$ defined by

\begin{equation}\label{e2}\mu^{(n)}(\sigma_n)=Z_n^{-1}\exp\left(-\beta H(\sigma_n)
+\sum_{x\in W_n}h_{\sigma(x),x}\right),\end{equation}\\

\begin{equation}\label{e3}Z_n=\int\!\!\!...\!\!\!\!\!\int\limits_{\Omega^{(p)}_{V_{n-1}}} \exp\left(-\beta
H({\widetilde\sigma}_n) +\sum_{x\in
W_{n}}h_{{\widetilde\sigma}(x),x}\right)
\lambda^{(p)}_{V_{n-1}}({d\widetilde\s_n}),\end{equation} where
$$\underbrace{\Omega_{W_{n}}\times\Omega_{W_{n}}\times...\times\Omega_{W_{n}}}_{3\cdot
2^{p-1}}=\Omega^{(p)}_{W_{n}},\ \ \
\underbrace{\lambda_{W_{n}}\times\lambda_{W_{n}}\times...\times\lambda_{W_{n}}}_{3\cdot
2^{p-1}}=\lambda^{(p)}_{W_{n}}, \ n,p\in \mathbb{N},$$ Let
$\sigma_{n-1}\in\Omega_{V_{n-1}}$ and
$\sigma_{n-1}\vee\omega_n\in\Omega_{V_n}$ is the concatenation of
$\sigma_{n-1}$ and $\omega_n.$ For $n\in \mathbb{N}$ we say that
the probability distributions $\mu^{(n)}$ are compatible if
$\mu^{(n)}$ satisfies the following condition:\vskip 0.1 truecm
\begin{equation}\label{e4}\int\!\!\!\!\!\!\!\!\!\int\limits_{\Omega_{W_n}\times\Omega_{W_n}}
\mu^{(n)}(\sigma_{n-1}\vee\omega_n)(\lambda_{W_n}\times
\lambda_{W_n})(d\omega_n)=
\mu^{(n-1)}(\sigma_{n-1}).\end{equation}\vskip 0.1 truecm

By Kolmogorov's extension theorem there exists a unique measure
$\mu$ on $\Omega_V$ such that, for any $n$ and
$\sigma_n\in\Omega_{V_n}$, $\mu \left(\left\{\sigma
|_{V_n}=\sigma_n\right\}\right)=\mu^{(n)}(\sigma_n)$. The measure
$\mu$ is called {\it splitting Gibbs measure} corresponding to
Hamiltonian (\ref{e1}) and function $x\mapsto h_x$, $x\neq x^0$.\\
Denote
\begin{equation}\label{e20}K(u,t,v)=\exp\left\{J_{3}\beta\xi_{1}\left(t,u,v\right)+J\beta\xi_{2}\left(u,v\right)
+J_{1}\beta\left(\xi_{3}\left(t,u\right)+\xi_{3}\left(t,v\right)\right)+\alpha\beta(u+v)\right\},\end{equation}
and
$$f(t,x)=\exp(h_{t,x}-h_{0,x}), \ \ (t,u,v)\in [0,1]^{3},\ x\in
V\setminus\{x^{0}\}.$$\vskip 0.3truecm

The following statement describes conditions on $h_x$ guaranteeing
compatibility of the corresponding distributions
$\mu^{(n)}(\sigma_n).$

 \begin{thm}\label{th1}  The measure
$\mu^{(n)}(\sigma_n)$, $n=1,2,\ldots$ satisfies the consistency
condition (\ref{e4}) iff for any $x\in V\setminus\{x^0\}$ the
following equation holds:
\begin{equation}\label{e5} f(t,x)=\prod_{>y,z<\in S(x)}
\frac{\int_0^1\int_0^1K(t,u,v)f(u,y)f(v,z)dudv}{\int_0^1\int_0^1K(0,u,v)f(u,y)f(v,z)dudv},
\end{equation} where $S(x)=\{y,z\},\ <y,x,z>$ is a ternary neighbor and $du=\lambda(du)$ is the Lebesgue measure
\end{thm}

\section{Lyapunov's operator $\mathcal L$ with degenerate kernel}

Now we consider the case $J_{3}\neq 0, \ J=J_{1}=\alpha=0$ for the
model (\ref{e1}) in the class of translational-invariant functions
$f(t,x)$ i.e $f(t,x)=f(t),$ for any $x\in V$. For such functions
equation (\ref{e1}) can be written as
\begin{equation}\label{e23}f(t)=\frac{\int_0^1\!\!\int_0^1K(t,u,v)f(u)f(v)dudv}{\int_0^1\!\!
\int_0^1K(0,u,v)f(u)f(v)dudv},
\end{equation}
 where $K(t,u,v)=\exp\left\{J_{3}\beta\xi_{1}\left(t,u,v\right)+J\beta\xi_{2}\left(u,v\right)
+J_{1}\beta\left(\xi_{3}\left(t,u\right)+\xi_{3}\left(t,v\right)\right)+\alpha\beta(u+v)\right\},$
$f(t)>0, \ t,u\in [0,1].$\\
We shall find positive continuous solutions to (\ref{e23}) i.e.
such that $f\in C^+[0,1]=\{f\in C[0,1]: f(x)\geq 0\}$.

 Define a nonlinear operator $H$ on the cone of positive continuous
functions on $[0,1]:$
$$
(Hf)(t)=\frac{\int_0^1\int_0^1K(t,s,u)f(s)f(u)dsdu}{\int_0^1\int_0^1
K(0,s,u)f(s)f(u)dsdu}.$$

 We'll study the existence of positive
fixed points for the nonlinear operator $H$ (i.e., solutions of
the equation (\ref{e23})). Put $ C_0^+[0,1]=C^+[0,1]\setminus
\{\theta\equiv 0\}.$ Then the set $C^+[0,1]$ is the cone of
positive continuous functions on $[0,1].$

We define the Lyapunov integral operator $\mathcal{L}$ on $C[0,1]$
by the equality (see \cite{newk})
$$\mathcal{L}f(t)=\int_0^1K(t,s,u)f(s)f(u)dsdu.$$
Put
$$\mathcal M_0=\left\{f\in C^+[0,1]: f(0)=1\right\}.$$
Denote by $N_{fix.p}(H)$ and $N_{fix.p}(\mathcal{L})$ are the set
of positive numbers of nontrivial positive fixed points of the
operators $N_{fix.p}(H)$ and $N_{fix.p}(\mathcal{L})$,
respectively.
\begin{thm}\label{l.1}\cite{newh}
\item{i)} The equation
\begin{equation}\label{e2.1}
Hf=f, \,\ f\in C^{+}_{0}[0,1]
\end{equation}
 has a positive solution iff the Lyapunov  equation
\begin{equation}\label{e2.2}
\mathcal{L}g=\lambda g, \,\ g\in C^{+}[0,1]
\end{equation}
has a positive solution in $\mathcal M_0$ for some $\lambda>0$.

\item{ii)} The equation  $Hf=f$ has a nontrivial positive solution iff
the Lyapunov equation $\mathcal{L}g=g$ has a nontrivial positive
solution.
\item{iii)} The equation
\begin{equation}\label{e2.3}
 \mathcal{L}f=\lambda f, \,\,\ \lambda>0
\end{equation}
 has at least one solution in $C_0^+[0,1].$
\item{iv)} The equation (\ref{e2.1}) has at least one solution in
$C_0^+[0,1]$. \item{v)} The equality
$N_{fix.p}(H)=N_{fix.p}(\mathcal{L})$ is hold.
\end{thm}

Let $\varphi_1(t), \,\ \varphi_2(t)$ and $\psi_1(t), \,\
\psi_2(t)$ are positive functions from $C_0^+[0,1]$.  We consider
Lyapunov's operator $\mathcal{L}$

\begin{equation}\label{l}(\mathcal{L}f)(t)=\int^{1}_{0}(\psi_{1}(t)\varphi_{1}(u)+\psi_{2}(t)\varphi_{2}(v))f(u)f(v)dudv.\end{equation}

and quadratic operator  $P$ on $\mathbb{R}^{2}$ by the rule

$$P(x,y) = (\alpha_{11}x^{2}+\alpha_{12}xy+\alpha_{22}y^{2}, \,\ \beta_{11}x^{2}+\beta_{12}xy+\beta_{22}y^{2}).$$

$$\alpha_{11}=\int_{0}^{1}\int_{0}^{1}\psi_{1}(u)\psi_{1}(v)\varphi_{2}(v)dudv,
\ \ \
\alpha_{12}=\int_{0}^{1}\int_{0}^{1}(\psi_{1}(v)\psi_{2}(u)+\psi_{1}(u)\psi_{2}(v))\varphi_{2}(v)dudv$$

$$\alpha_{22}=\int_{0}^{1}\int_{0}^{1}\psi_{2}(u)\psi_{2}(v)\varphi_{2}(v)dudv,
\ \ \
\beta_{11}=\int_{0}^{1}\int_{0}^{1}\psi_{1}(u)\psi_{1}(v)\varphi_{1}(u)dudv,$$

$$\beta_{12}=\int_{0}^{1}\int_{0}^{1}(\psi_{1}(u)\psi_{2}(v)+\psi_{1}(v)\psi_{2}(u))\varphi_{1}(u)dudv, \ \ \
\beta_{22}=\int_{0}^{1}\int_{0}^{1}\psi_{2}(u)\psi_{2}(v)\varphi_{1}(u)dudv.$$

\begin{lemma}\label{l3.1.} The Lyapunov's operator
$\mathcal{L}$ has a nontrivial positive fixed point iff the
quadratic operator $P$ has a nontrivial positive fixed point,
moreover $N_{fix}^{+}(H_k)=N_{fix}^{+}(P)$.
\end{lemma}

\proof $a)$ Put
$$\mathbb{R}_{2}^{+}=\{(x,y)\in \mathbb{R}^{2}: \,\ x\geq0, y\geq0\}, \ \mathbb{R}_{2}^{>}=\{(x,y)\in \mathbb{R}^{2}: \,\ x>0,
y>0\}.$$ Let $f(t)\in C_0^+[0,1]$ be a nontrivial positive fixed
point of $\mathcal{L}$. Let
$$c_{1}=\int^{1}_{0}\varphi_1(u)f(u)f(v)dudv, \ \
c_{2}=\int^{1}_{0}\varphi_2(u)f(u)f(v)dudv$$ Clearly, $c_1>0, \,\
c_2>0$ and $f(t)=c_1\psi_1(t)+c_2\psi_2(t)$. If we put
$f(t)=c_1\psi_1(t)+c_2\psi_2(t)$ to the equation (\ref{l}) we'll
get
$$c_1=\alpha_{11}c_1^{2}+\alpha_{12}c_1c_2+\alpha_{22}c_2^{2}, \ c_2=\beta_{11}c_1^{2}+\beta_{12}c_1c_2+\beta_{22}c_2^{2}.$$
Therefore, the point $(c_1, c_2)$ is fixed point of the quadratic
operator $P.$

$b)$ Assume, that the point $(x_0, y_0)$ is a nontrivial positive
fixed point of the quadratic operator $P,$ i.e. $(x_0, y_0)\in
\mathbb{R}_{2}^{+}\setminus\{\theta\}$ and numbers $x_0, y_0$
satisfies following equalities

$$\alpha_{11}x_0^{2}+\alpha_{12}x_0y_0+\alpha_{22}y_0^{2}=x_0, \ \beta_{11}x_0^{2}+\beta_{12}x_0y_0+\beta_{22}y_0^{2}=y_0.$$

Similarly, we can prove that the function
$f_0(t)=x_0\psi_1(t)+y_0\psi_2(t)$ is a fixed point of the
operator $\mathcal{L}$ and $f_0(t)\in C_0^+[0,1]$. This completes
the proof.
\endproof

We define positive quadratic operator $\mathcal{Q}$:
$$ \mathcal{Q}(x,y) = (a_{11}x^{2}+a_{12}xy+a_{22}y^{2}, \,\ b_{11}x^{2}+b_{12}xy+b_{22}y^{2}).$$
\begin{pro}\label{l4.1.} \item{i)} If $\omega=(x_0,y_0)\in \mathbb{R}_{2}^{+}$ is a positive fixed point of $\mathcal{Q}$,
then $\lambda_0=\frac{x_0}{y_0}$ is a root of the following
equation
\begin{equation}\label{e4.1}
a_{11}\lambda^{3}+(a_{12}-b_{11})\lambda^{2}+(a_{22}-b_{12})\lambda-b_{22}=0.
\end{equation}
\item{ii)}  If the positive number $\lambda_0$
is a positive root of the equation (\ref{e4.1}), then the point
$\omega_0=(\lambda_{0}y_0,y_0)$ is a positive fixed point of
$\mathcal{Q}$, where
$y^{-1}_{0}=a_{11}+a_{12}\lambda_{0}+a_{22}\lambda_{0}^{2}.$\end{pro}

\proof {\it i)} Let the point $\omega=(y_0,x_0)\in
\mathbb{R}_{2}^{+}$ be a fixed point of $\mathcal{Q}$. Then $$
a_{11}x_{0}^{2}+a_{12}x_{0}y_{0}+a_{22}y_{0}^{2}=x_{0}, \,\,\
b_{11}x_{0}^{2}+b_{12}x_{0}y_{0}+b_{22}y_{0}^{2}=y_{0}$$  Using
the equality $\frac{x_{0}}{y_0} = \lambda_0$ we obtain
$$a_{11}\lambda_{0}^{2}y_{0}^{2}+a_{12}\lambda_{0}y_{0}^{2}+a_{22}y_{0}^{2}=\lambda_{0}y_{0},
\,\,\
b_{11}\lambda_{0}^{2}y_{0}^{2}+b_{12}\lambda_{0}^{2}y_{0}^{2}+b_{22}y_{0}^{2}=y_{0}.$$
Thus we get
$$ \frac{a_{11}\lambda_{0}^{2}+a_{12}\lambda_{0}+a_{22}}{b_{11}\lambda_{0}^{2}+b_{12}\lambda_{0}+b_{22}}=\lambda_{0}.$$
Consequently,
$$ a_{22}+(a_{12}-b_{22})\lambda_{0}+(a_{11}-b_{12})\lambda_{0}^{2}-b_{11}\lambda_{0}^{3}=0. $$
{\it ii)} Let $\lambda_0>0$ is a root of the cubic equation
(\ref{e4.1}). Put $x_0=\lambda_0y_0$, where
$$x_0=\frac{\lambda_0}{a_{11}\lambda_{0}^{2}+2a_{12}\lambda_{0}+a_{22}}.$$
Since
$$a_{11}x_{0}^{2}+2a_{12}x_{0}y_{0}+a_{22}y_{0}^{2}=\frac{1}{a_{11}\lambda_{0}^{2}+2a_{12}\lambda_{0}+a_{22}},$$
we get
$$a_{11}x_{0}^{2}+2a_{12}x_{0}y_{0}+a_{22}y_{0}^{2}=y_0.$$
By the other hand
$$ a_{22}+(a_{12}-b_{22})\lambda_{0}+(a_{11}-b_{12})\lambda_{0}^{2}-b_{11}\lambda_{0}^{3}=0. $$
Then we get
$$b_{11}\lambda_{0}^{2}+b_{12}\lambda_{0}+b_{22}=
\lambda_{0}(a_{11}\lambda_{0}^{2}+a_{12}\lambda_{0}+a_{22}).$$
From the last equality we get
$$\frac{\lambda_{0}}{a_{11}\lambda_{0}^{2}+a_{12}\lambda_{0}+a_{22}}=\frac{b_{11}\lambda_{0}^{2}+b_{12}
\lambda_{0}+b_{22}}{(a_{11}\lambda_{0}^{2}+a_{12}\lambda_{0}+a_{22})^{2}}=$$
$$=b_{11}x_{0}^{2}+2b_{12}x_{0}y_{0}+b_{22}y_{0}^{2}=y_{0}.$$
This completes the proof.
\endproof
Denote
$$P(\lambda)=\alpha_{11}\lambda^{3}+(\alpha_{12}-\beta_{11})\lambda^{2}+(\alpha_{22}-\beta_{12})\lambda-\beta_{22}=0, \
 \mu_0=\alpha_{11}, \,\ \mu_1=\alpha_{12}-\beta_{11}, \,\
\mu_2=\alpha_{22}-\beta_{12}, \,\ \mu_3=\beta_{22},$$
\begin{equation}\label{add} P_3(\xi)=\mu_0\xi^{3}+\mu_1\xi^{2}+\mu_2\xi-\mu_3,\end{equation}
$$D=\mu_{1}^{2}-3\mu_0\mu_2, \
\alpha=-\frac{\mu_1+\sqrt{D}}{3\mu_0}, \,\,\
\beta=-\frac{\mu_1-\sqrt{D}}{3\mu_0}.$$

\begin{thm}\label{t4.1.} Let $\mathcal{Q}$ satisfies one of the following
conditions
\item{i)} $D\leq0;$
\item{ii)} $D>0, \beta\leq0;$
\item{iii)} $D>0, \alpha\leq0, \beta>0;$
\item{iv)} $D>0, \alpha>0, P_3(\alpha)<0;$
\item{v)} $D>0, \alpha>0, P_3(\alpha)>0, P_3(\beta)>0,$
 then $\mathcal{Q}$ has a unique nontrivial positive fixed
point.
\end{thm}

\proof The proof of Theorem \ref{t4.1.} is basis on monotonous
property of the function $P_3(\xi)$. Clearly,

\begin{equation}\label{e4.3}
\left(P_3(\xi)\right)^{'}=3\mu_0\xi^{2}+2\mu_1\xi+\mu_2.
\end{equation}
and $P'_{3}(\alpha)=P'_{3}(\beta)=0.$ Moreover,

 $i)$ In the case $D\leq0,$  by the equality (\ref{e4.3}) the
function $P_3(\xi)$ is an increasing function on $\mathbb{R}$ and
$P_3(0)=-b_{11}<0.$ Therefore, the polynomial $P_3(\xi)$ has a
unique positive root.

$ii)$ Let $D>0$ and $\beta\leq0.$ For the case $D>0$ the function
$P_3(\xi)$ is an increasing function on $(-\infty,
\alpha)\cup(\beta, \infty)$ and it is a decreasing function on
$(\alpha, \beta)$. Hence from the inequality $P_3(0)<0$ the
polynomial $P_3(\xi)$ has a unique positive root.

$iii)$ Let $D>0, \alpha\leq0$ and $\beta>0.$ Since  the function
$P_3(\xi)$ is decreasing on $(\alpha, \beta)$ and increasing on
$(\beta, \infty)$ the polynomial $P_3(\xi)$ has a unique positive
root as $P_3(0)<0.$

$iv)$ Let $D>0, \alpha>0$ and $P_3(\alpha)<0.$ Then $\max_{\xi\in
(-\infty,\beta)}P_3(\xi)=P_3(\alpha)<0$. Consequently, by the
function $P_3(\xi)$ is increasing on $(\beta, \infty)$ the
polynomial $P_3(\xi)$ has a unique positive root $\xi_0 \in
(\beta, \infty)$.

$v)$ Let $D>0, \alpha>0, P_3(\alpha)>0$ and $P_3(\beta)>0.$ Then
$\min_{\xi\in (\alpha,\infty)}P_3(\xi)=P_3(\beta)>0$. From the
function $P_3(\xi)$ on $(-\infty, \alpha)$,  $P_3(\xi)$ $P_3(\xi)$
has a unique positive root $\xi_0 \in (0, \alpha)$, as $P_3(0)<0$
and $P_3(\alpha)>0$.

From the upper analysis  and by Lemmas \ref{l4.1.} it follows that
the Theorem  \ref{t4.1.}.
\endproof

\begin{thm}\label{t4.2.} Let be $ D>0$. If $\mathcal{Q}$ satisfies one of the following
conditions
\item{i)} $\alpha>0, P_3(\alpha)=0, P_3(\beta)<0;$
\item{ii)} $\alpha>0, P_3(\alpha)>0, P_3(\beta)=0,$
then QO $\mathcal{Q}$ has two nontrivial positive fixed points and
$N_{fix}^{+}(\mathcal{Q})=N_{fix}^{>}(\mathcal{Q})=2$.
\end{thm}

\proof

$i)$ Let $\alpha>0, P_3(\alpha)=0$ and $P_3(\beta)<0.$ Then
$\max_{\xi\in (-\infty,\beta)}P_3(\xi)=P_3(\alpha)=0$ and
$\xi_1=\alpha$ is the root of the polynomial $P_3(\xi)$. By the
increase  property on $(\beta, \infty)$ of the function $P_3(\xi)$
the polynomial $P_3(\xi)$ has a  root $\xi_2 \in (\beta, \infty)$,
as $\beta>0$ and $P_3(\beta)<0$. There is not any other positive
roots of the polynomial $P_3(\xi)$.\\
 $ii)$ Let $\alpha>0,
P_3(\alpha)>0$ and $P_3(\beta)=0.$ Then by the increase property
on $(-\infty, \alpha)$ of the function $P_3(\xi)$ the polynomial
$P_3(\xi)$ has a  root $\xi_1 \in (0, \alpha)$. By the other hand
$\min_{\xi\in(\alpha,\infty)}P_3(\xi)=P_3(\beta)=0$ and the number
$\xi_2=\alpha$ is the second positive root of the polynomial
$P_3(\xi)$. The polynomial $P_3(\xi)$ has not another roots. From
above  and by Lemmas \ref{l4.1.} we get Theorem \ref{t4.2.}.
\endproof

\begin{thm}\label{t4.2.} Let be $ D>0$. If $\mathcal{Q}$ satisfies one of the following
conditions
\item{i)}  $\alpha>0, P_3(\alpha)=0, P_3(\beta)<0;$
\item{ii)} $\alpha>0, P_3(\alpha)>0, P_3(\beta)=0,$\\
then $\mathcal{Q}$ has two nontrivial positive fixed points and
$N_{fix}^{+}(\mathcal{Q})=N_{fix}^{>}(\mathcal{Q})=2$.
\end{thm}
\proof $i)$ Let $\alpha>0, P_3(\alpha)=0$ and $P_3(\beta)<0.$ Then
$\max_{\xi\in (-\infty,\beta)}P_3(\xi)=P_3(\alpha)=0$ and
$\xi_1=\alpha$ is the root of the polynomial $P_3(\xi)$. By the
increase  property on $(\beta, \infty)$ of the function $P_3(\xi)$
the polynomial $P_3(\xi)$ has a  root $\xi_2 \in (\beta, \infty)$,
as $\beta>0$ and $P_3(\beta)<0$. There is not any other positive
roots of the polynomial $P_3(\xi)$.\\
 $ii)$ Let $\alpha>0,
P_3(\alpha)>0$ and $P_3(\beta)=0.$ Then by the increase property
on $(-\infty, \alpha)$ of the function $P_3(\xi)$ the polynomial
$P_3(\xi)$ has a  root $\xi_1 \in (0, \alpha)$. By the other hand
$\min_{\xi\in(\alpha,\infty)}P_3(\xi)=P_3(\beta)=0$ and the number
$\xi_2=\alpha$ is the second positive root of the polynomial
$P_3(\xi)$. The polynomial $P_3(\xi)$ has not another roots. From
above  and by Lemmas \ref{l4.1.} we get Theorem \ref{t4.2.}.
\endproof

\end{document}